\documentclass[12pt]{amsart}
\theoremstyle{plain}
\newtheorem{theorem}{Theorem}[section]
\newtheorem{corollary}[theorem]{Corollary}
\newtheorem{lemma}[theorem]{Lemma}
\newtheorem{proposition}[theorem]{Proposition}

\theoremstyle{definition}
\newtheorem{definition}[theorem]{Definition}

\theoremstyle{remark}
\newtheorem{remark}[theorem]{Remark}

\usepackage{amssymb}
\def\ker{\operatorname{ker}}
\def\Ind{\operatorname{Ind}}
\def\deg{\operatorname{deg}}
\def\EXT{\operatorname{\it{Ext}}}
\def \Ext #1 #2{\EXT(#1,#2)}
\def \KK {\operatorname{KK}}
\def \Hom {\operatorname{Hom}}
\def \ot {\otimes}
\def \mb {\mathbb}
\def \mc {\mathcal}
\begin{document}
\title [ A note on Kasparov product and Duality ]
       {A note on Kasparov product and Duality}

\author { Hyun Ho  \quad Lee }
\address {Department of Mathematical Sciences
         Seoul National University
         Seoul, South Korea 151-747 }
\email{hyun.lee7425@gmail.com}
\date{September 5, 2010}
 \begin{abstract}

  Using Paschke-Higson duality, we can get a natural index pairing $K_{i}(A) \times K_{i+1}(D_{\Phi}) \to
  \boldsymbol{Z} \quad  (i=0,1) (\mbox{mod}2)$, where $A$ is a separable
  $C\sp*$-algebra, and $\Phi$ is a representation of $A$ on a separable infinite dimensional Hilbert space $H$.
  It is proved that this is a special case of the Kasparov Product. As a step, we show a proof of Bott-periodicity for KK-theory
  asserting that $\mathbb{C}_1$ and $S$ are $KK$-equivalent using the odd index pairing.

\end{abstract}
\maketitle

\section{Introduction} \label{S:intro}
In \cite{hi1}, Higson showed a connection between K-theory of the
essential commutant of a $C\sp{*}$-algebra $A$ and the extension
group of A by the $C\sp{*}$-algebra of the compact operators on a
separable infinite dimensional Hilbert space and defined  index pairings (dualities)
\begin{equation}\label{E:duality}
K_{i}(A)\times K_{i+1}(D_{\Phi}) \to
  \mathbb{Z}
\end{equation}  for $i=0,1 (\mbox{mod}2)$ based on an earlier work of Paschke \cite{pa}.

In light of this duality, K-theory of the essential commutant of a $C\sp*$-algebra $A$ can be viewed as K-homology of $A$ (See Theorem 1.5 in \cite{hi1}). Since Kasparov's KK-theory generalizes both K-theory and K-homology, we can expect that each index pairing given by (\ref{E:duality}) is a Kasparov product.

It is a goal of our paper to show that Paschke-Higson duality can be realized as a Kasparov product.  For this we need to use Bott periodicity for $\KK$-theory of the form $\KK(A,B)\cong \KK(SA,SB)$ so that we give a proof of Bott periodicity using the odd pairing as an application.
\section{Paschke-Higson
Duality and Index pairing} \label{S:Pairing}
In this section, we review the Paschke-Higson duality theory \cite{hi1} more carefully
for the convenience of the reader.\\
Throughout this article we denote by $H$ a separable
infinite dimensional Hilbert space, by $B(H)$ the
set of linear bounded operators on $H$, by
$K(H)$(or just $K$) the ideal of
compact operators on $H$, and
by $\mathcal{Q}(H)$(or just $\mathcal{Q}$) the Calkin algebra.\\
We use the following notation : if $X$ and $Y$ are operators in
$B(H)$ we shall write
$$ X \sim Y $$ if $X$ and $Y$ differ by a compact operator. \\

 Note that every *-representation $\Phi$ of $A$ on $H$ induces a
*-homomorphism $\dot{\Phi}$ of $A$ into the Calkin algebra.
\begin{definition}\label{D:admissible}
Let $A$ be a $C^{*}$-algebra. A *-representation $\Phi: A \to
B(H)$ is called \emph{admissible} if it is non-degenerate
and $\ker(\dot{\Phi})=0$.
\end{definition}
\begin{remark}
If a *-representation is admissible, then it is necessarily faithful
and its image contains no non-zero compact operators.
\end{remark}
\begin{definition}
Let $\Phi$ be a *-representation of $A$ on $H$. We define
the \emph{essential commutant} of $\Phi(A)$ in
$B(H)$ as
$$ D_{\Phi}(A)=\{T\in B(H) \mid [\Phi(a),T]\sim 0 \quad \text{for all}\quad a \in A\}$$
\end{definition}
Given two representations $\Phi_{0}$ and $\Phi_{1}$ on
$H_{0}$ and $H_{1} $ respectively, we say they are
\emph{approximately unitarily equivalent} if there exists a sequence
$\{U_{n} \}$ consisting of unitaries in
$B(H_{0},H_{1})$ such that for any $a \in A$
\begin{align*}
U_{n}\Phi_{0}(a)U_n^* &\sim \Phi_{1}(a),  \\
\|U_{n}\Phi_{0}(a)U_n^* &- \Phi_{1}(a)\| \rightarrow 0 \quad \mbox{as} \quad n \to \infty.
\end{align*}
We write $\Phi_{0} \sim_{u} \Phi_{1}$ in this case.
\begin{theorem}(Voiculescu)
Let $A$ be a separable $C^{*}$-algebra and  $\Phi_{i}:A \to
B(H_{i})$  $(i=0,1) $ be non-degenerate
*-representations. Then if $\ker \dot{\Phi}_{0}\subset \ker
\Phi_{1}$, then $\Phi_{0}\oplus\Phi_{1} \sim_{u}\Phi_{0}$.
\end{theorem}
 \begin{proof}
  See Corollary 1 in p343 of \cite{ar}.
 \end{proof}
 \begin{corollary}\label{C:equivalence1}
 Assume $\Phi_{i}:A\to B(H_{i}) $ are
 admissible representations for $i=0,1$. Then
 $\Phi_{0}\sim_{u}\Phi_{1}$.
 \end{corollary}
 \begin{proof}
 Admissibility implies $\ker\dot{\Phi}_{i}= 0$. By symmetry, the result follows.
 \end{proof}
  Recall that an extension of a unital separable $C\sp{*}$-algebra $A$ is a unital
  *-monomorphism $\tau$ of $A$ into the Calkin algebra. We say $\tau$ is \emph{split}
  if there is a non-degenerate *-representation $\rho$ such that
  $\dot{\rho}=\tau$ and \emph{semi-split} if there is a completely positive map $\rho$
  such that $ \dot{\rho}=\tau$. In general, an extension $\tau$ of $A$ by $B$ is a $*$-homomorphism from $A$ to $Q(B)$ where $Q(B)=M(B)/B$ and $M(B)$ is the multiplier algebra of $B$.  We say  $\tau$ is \emph{split} or \emph{semi-split} if the lifting $\rho$ is a *-homomorphism to $M(B)$ or a completely positive map to $M(B)$ respectively. When $B$ is stable, $\EXT(A,B)$ is the quotient of $\Hom(A,Q(B))$ by the equivalence relation generated by addition of trivial elements and unitary equivalence.  
 \begin{corollary}\label{C:equivalence2}
 Let $A$ be a separable unital $C^{*}$-algebra. If $\tau$ is a unital injective
 extension of $A$ and if $\sigma$ is a split unital extension of $A$,
 then $\tau\oplus\sigma$ is unitarily equivalent to $\tau$.
 \end{corollary}
 \begin{proof}
 See p352-353 in \cite{ar}.
 \end{proof}
  Now we will prove that there is at least one admissible representation of $A$.
  \begin{proposition}
   There is a non-degenerate *-representation $\Phi$ of for a separable
   $C\sp{*}$-algebra $A$ such that $\ker \dot{\Phi}=0$.
  \end{proposition}
  \begin{proof}
   Let $\pi$ be a faithful representation of $A$ on $\mathbb{H_{\pi}}$. Take
   $\Phi$ to be $\pi^{(\infty)}=\pi \oplus \pi \oplus \cdots$ and
   $H=H_{\pi}^{(\infty)}=H_{\pi} \oplus
   H_{\pi} \oplus \cdots$.
  \end{proof}
  \begin{proposition}\label{P:isodual}
  If $\Phi$, $\Psi$ are admissible representations of $A$ on
  $H$,  $D_{\Phi}(A)$  is isomorphic to $D_{\Psi}(A)$.
  \end{proposition}
   \begin{proof}
    By Corollary \ref
  {C:equivalence1}, there is a unitary $U \in B(H)$ such that for any $a \in A$
  \[U\Phi(a)U^* \sim \Psi(a).\]
  We define a map $\theta$ on $D_{\Phi}(A)$ by
  $\theta(T)=UTU^*$ for $T \in D_{\Phi}(A)$. We note that for any $a \in A$
  \begin{align*}
 \Psi(a)\theta(T)- \theta(T)\Psi(a)&=\Psi(a) UTU^*-UTU^*\Psi(a)\\
                                   &\sim U\Phi TU^*- U T \Phi U^*\\
                                   &\sim U (\Phi T-T\Phi)U^*\\
                                   &\sim 0.
  \end{align*}
  Thus $\theta(T) \in D_{\Psi}(A)$ and $\theta$ is an isomorphism onto $D_{\Psi}(A)$ since $\theta=Ad (U)$.

  \end{proof}
 \begin{definition}\label{D:Dualalgebra}
  When $\Phi$, $\Psi$ are \emph{admissible} representations of $A$ on
  $H$,  $D_{\Phi}(A)$  is isomorphic to $D_{\Psi}(A)$ by  Proposition \ref
  {P:isodual}. Thus
  we define $D(A)=D_{\Phi}(A)$ as the dual algebra of
  $A$ up to unitary equivalence.
 \end{definition}
 If $p$ is a projection in $D_{\Phi}(A)$, we call it \emph{ample} and can define an
 extension
 $$\tau=\tau_{\Phi,p}: A \xrightarrow{p \Phi(\bullet ) p = \Phi_{p}} B(pH)\xrightarrow{\pi}\mathcal{Q}(pH)$$
  In general, when $A$ is separable and $B$ is stable, we can express invertible elements $\tau$ in $\EXT(A,B)$ as pairs $(\phi,p)$, where $\phi$ is a representation from $A$ to $M(B)$ and $p$ is a projection in $M(B)$ which commutes with $\phi(A)$ modulo $B$, i.e, $\tau(\cdot)=q_B(p\phi(\cdot)p)$ where $q_B= M(B) \to Q(B) $. We denote $\EXT^{-1}(A,B)$ by the group of invertible elements in $\EXT(A,B)$ \cite{bl}, \cite{kas80}.

 To define Paschke-Higson duality, we need the
 following two technical lemmas.
 \begin{lemma}\label{L:K_0surjectivity}
 Let $A$ be a unital $C^{*}$-algebra. For any $ \alpha \in K_{0}(D(A))$, there
 exists an ample projection $p \in D(A)$ such that
 $\alpha=[p]_{0}$.
 \end{lemma}
 \begin{proof}
Step1: By Corollary \ref{C:equivalence2}, there is a unitary $u \in
B(pH\oplus H,H)$ such that Ad
$(u) (\Phi_{p}\oplus \Phi) (a) \sim \Phi(a)\, \quad \mbox{for any}
\quad a \in A$ if $p$ is ample. Let $ U = \binom{p}{0} u $. We can
easily check that $U \in \mathbb{M}_{2}(D_{\Phi}(A))$ and $UU^{*}=
\left(
\begin{matrix}
p & 0 \\
0 & 0
 \end{matrix}
\right), U^{*}U=
 \left(
\begin{matrix}
p & 0 \\
0 & I
\end{matrix}
\right)$. Therefore we have $[p \oplus I]_{0}=[p\oplus 0]_{0}$. This implies $[p]_{0}+ [I]_{0}=[p]_{0}$.
In particular, $[I]_{0}=0 $. Similarly, we can conclude every
two ample projections are Murray-von Neumann equivalent. \\
 Step2: Note that $p\oplus 1$ is always ample whether $p$ is ample or not
 because $(\Phi \oplus \Phi)_{p\oplus 1}(a)$ is never compact unless $a$ is zero in $A$.\\
 Step3: Any element in $K_{0}(D_{\Phi}(A))$ can be written as $[q]_{0}-[I_{n}]_{0}$ for some $q \in \mathbb{M}_{n}(D_{\Phi}(A))$.
 As we observed in Step1, this is just $[q]_{0}$. By Step2, we may assume $q$ is ample for
 $\Phi^{(n)}=\underbrace{\Phi\oplus\Phi\oplus \cdots \oplus \Phi}_{n}$ .  Now if we can show $[q]_{0}=[p]_{0}$ for some \emph{ample} $p \in D_{\Phi}(A)$, we are done.\\
 Since  $\Phi^{(n)} \sim_{u} \Phi $, there exists $v: H^{n} \to H$
 such that
 \begin{enumerate}
\item $v^{*}v=1_{B(H^{n})}, vv^{*}=1_{B(H)}$ \\
\item Ad$(v)\Phi^{(n)}(a)-\Phi(a)\in K \quad  \text{for any} \quad a \in A $.
\end{enumerate}
Then $[q]_{0}=[vqv^{*}]_{0}$. In addition, $[vqv^{*},\Phi(a)]\sim
v[q,\Phi^{(n)}]v^* \sim 0 $ for every $a \in A$. Thus $vqv^*$ is
ample and we are done.
\end{proof}

 \begin{lemma}\label{L:K_1surjectivity}
 Let $A$ be as above. For any $ \alpha \in K_{1}(D(A))$, there
 exists an unitary $u \in D(A)$ such that
 $\alpha=[u]_{1}$.
 \end{lemma}
 \begin{proof}
 Assume $U \in \mathbb{M}_{n}(D_{\Phi}(A)) \approx D_{\phi^{(n)}}(A)$ is a unitary which represents $\alpha \in K_{1}(D_{\Phi}(A))$.
 Let $V$ be $(\overbrace{1,0, \cdots,0}^{n})^{T}v$ where $v$ is defined in Lemma \ref{L:K_0surjectivity}
 and $S=\left(\begin{matrix}
V& 1-VV^{*} \\
 0 & V^{*}   \end{matrix} \right) \in \mathbb{M}_{2n}(D_{\Phi}(A)) $. It is easy to check that
 $VUV^{*}+1-VV^{*}=\left(\begin{matrix}
vUv^{*}& 0 \\
 0     & 1 \end{matrix} \right)$ and $S\left(\begin{matrix}
U & 0 \\
 0     & 1 \end{matrix} \right)S^{*}=\left(\begin{matrix}
VUV^{*}+1-VV^{*}& 0 \\
 0     & 1 \end{matrix} \right)$. Therefore $[U]_{1}=[vUv^{*}]_{1}$ and  $[vUv^{*},\Phi(a)] \sim v[U,
 \Phi^{n}(a)]v^*\sim 0$
 for every $a \in A$. Thus we take $u$ as $vUv^*$.
\end{proof}

 \begin{remark}
  A unital $C\sp{*}$-algebra $A$ is said to have $K_1$-surjectivity if the
  natural map from $U(A)/U_{0}(A)$ to $K_{1}(A)$ is surjective and is
  said to have (strong) $K_0$-surjectivity if the group $K_{0}(A)$ is generated
  by $\{ [p]\mid p \,\,\mbox{is a projection in}\,\, A\}$. Thus Lemma
  \ref{L:K_0surjectivity} and Lemma \ref{L:K_1surjectivity} show $D_{\Phi}(A)$ has
   (strong) $K_0$-surjectivity and $K_1$-surjectivity.
 \end{remark}
 
 Now we are ready to define an index pairing between $K_{i}(A)$ and $K_{i+1}(D_{\Phi}(A))$ for $i=0,1$.
 In the following two definitions, $\Ind$ will denote the classical \emph{Fredholm
 index}. 

 Given a projection $ p \in \mathbb{M}_{k}(A)$ and a unitary $u \in D_{\Phi}(A)$,
 when $\Phi_{k}$ is k-th amplification of $\Phi$, the operator
 $$\Phi_{k}(p)u^{(k)}\Phi_{k}(p): \Phi_k(p)(H^k) \to
 \Phi_k(p)(H^k)$$
is essentially a unitary, and therefore \emph{Fredholm}.
\begin{definition}\label{D:EvenIndex}
The (even) index pairing $K_{0}(A) \times  K_{1}(D_{\Phi}(A)) \to
\qquad\mathbb{Z}$ is defined by
$$([p], [u])  \xrightarrow \qquad \Ind \left(\Phi_{k}(p)u^{(k)}
\Phi_{k}(p) \right)$$ where $p \in \mathbb{M}_{k}(A)$ and $\Phi_{k}$ is k-th amplification of $\Phi$.
\end{definition}
 Similarly, given a unitary $v \in \mathbb{M}_{k}(A)$ and a projection $p \in D_{\Phi}(A)$, the
 operator $$p^{(k)}\Phi_{k}(v)p^{(k)}
-(1-p^{(k)}): H^{k} \to H^{k} $$
is essentially a unitary, and therefore \emph{Fredholm}.
\begin{definition}\label{D:OddIndex}
 The (odd) index pairing $K_{1}(A) \times  K_{0}(D_{\Phi}(A)) \to
 \mathbb{Z}$ is defined by
$$([v], [p])  \xrightarrow \qquad \Ind \left(p^{(k)}\Phi_{k}(v)p^{(k)}
-(1-p^{(k)})\right)$$ where $ v \in \mathbb{M}_{k}(A)$ and $\Phi_{k}$ is k-th amplification of $\Phi$.
\end{definition}

\section{Kasparov Product and Duality} \label{S:Kasparov}
In this section, we prove our main results: Each index pairing is a
special case of the Kasparov product. Before doing this, we recall
some rudiments of KK-theory. We only give a brief description of elements of the Kasparov group $\KK(A,B)$, but for the complete information we refer the reader to the original source \cite{kas80} or to \cite{jeth} and \cite{bl} for detailed expositions.

 $\KK(A,B)$ is described in terms of  triples $(E,\phi,F)$, which we call cycles, where $E$ is a countably generated graded Hilbert B-module, $\phi: A \to \mc{L}_B(E)$ is a representation and $F \in \mc{L}_B(E)$ is an element of degree 1 such that $[F,\phi(a)],(F^2-1)\phi(a),(F^*-F)\phi(a)$ are in $\mc{K}_B(E)$ for any $a \in A$. We remind the reader that the commutator $[\, ,\, ]$ is graded and $\phi$ preserves the grading. We shall denote $\mb{E}(A,B)$ by the set of all cycles. A degenerate cycle is a triple $(E,\phi, F)$ such that $[F,\phi(a)]=(F^2-1)\phi(a)=(F^*-F)\phi(a)=0$ for all $a \in A$. An operator homotopy through cycles is a homotopy $(E,\phi, F_t)$, where $t \to F_t$ is norm continuous. A theorem of Kasparov \cite{kas80} shows that $\KK(A,B)$ is the quotient of $\mb{E}(A,B)$ by the equivalence relation generated by addition of degenerate cycles, unitary equivalence and operator homotopy.

A graded (maximal) tensor product of two graded $C\sp*$-algebras $C$ and $D$ denoted by $C\hat{\ot}D$ is the universal enveloping $C\sp*$-algebra of the algebraic tensor product $C \odot D$ with a new product and involution on $C \odot D$ by
\begin{align*}
(c_1\hat{\ot}d_1)(c_2\hat{\ot}d_2)&=(-1)^{\deg d_1 \deg c_2}(c_1c_2\hat{\ot}d_1d_2)\\
(c_1\hat{\ot}d_1)^*&=(-1)^{\deg c_1 \deg d_1}(c_1^*\hat{\ot}d_1^*)
\end{align*}
Let $\mc{E}=(E,\phi,F) \in \mb{E}(A,B)$. Then we can form a Hilbert $B \hat{\ot} C$-module $E \hat{\ot} C$. In addition, $\phi \hat{\ot} 1 : A\hat{\ot}C \to \mc{L}_{B\hat{\ot}C}(E\hat{\ot}C)$ is defined by \[\phi\hat{\ot}1(a\ot c)=\phi(a) \hat{\ot} c. \] It follows that $(E\hat{\ot}C, \phi\hat{\ot}1, F\hat{\ot}1) \in \mb{E}(A \hat{\ot} C, B\hat{\ot}C)$. We denote this map by $\tau_C: \KK(A,B) \to \KK(A \hat{\ot} C, B \hat{\ot} C)$. Thus the image of $\mc{E}$ is $\tau_C(\mc{E})$. 

 For $\mc{E}\in \mb{E}(A,B)$ and $\mc{F}\in \mb{E}(B,C)$, there exists a Kasparov product of $\mc{E}$ by $\mc{F}$, which will be denoted by $\mc{E} \cdot\mc{F}$, in $\mb{E}(A,C)$ \cite{kas80}. This product is unique up to operator homotopy so that we define the product $\KK(A,B)\times \KK(B,C) \to \KK(A,C)$ using the same notation.

For the remaining, $A$ and $B$ will denote (ungraded) separable $C\sp*$-algebras.
 \begin{proposition} \label{P:1}
 $KK(S,B)=K_{1}(B)$  where $S=\{f\in C(\mathbb{T})\mid f(1)=0\}$.
 \end{proposition}
 \begin{proof}
  It is well-known. We just note that any unitary in $K_{1}(B)$ can be lifted to $\boldsymbol{\phi} \in KK(S,B) $
  where $\phi:S \to B$ is determined by sending $z-1$ to $u-1$.
\end{proof}
 \begin{proposition}\label{P:KKCC}
 $\KK(\mb{C},\mb{C}) \cong \mb{Z}$.
 \end{proposition}
 \begin{proof}
 Every element in $\KK(\mb{C},\mb{C})$ can be represented by a triple
 \[ \left( \hat{H},\left( \begin{matrix}\phi & 0 \\
 0& \phi \end{matrix}\right), \left( \begin{matrix}0 & u^{*} \\
 u& 0 \end{matrix}\right)  \right)\] which satisfies  $[u,p] \sim 0 $, $(u-u^*)p \sim 0$, and $(u^2-1)p \sim 0$ for $\phi(1)=p$. Here $\hat{H}$ means $H\oplus H$ with the standard even grading. Such a triple is called a standard triple. Then a standard triple is mapped to $\Ind(pup)=\Ind(\phi(1)u\phi(1))$. And this map is a group isomorphism ( See Example 17.3.4 in
\cite{bl}).
 \end{proof}
Recall that $\mathbb{C}_{1}$ is the Clifford algebra of
$\mathbb{R}$,  which is isomorphic to $\mathbb{C}\oplus\mathbb{C}$
with the odd grading that transposes the two copies of $\mathbb{C}$. Thus $B\hat{\ot}\mb{C}_1$ is the $C\sp*$-algebra $B\oplus B$ with the odd grading.
\begin{definition}
\[\KK^1(A,B)=KK(A,B\hat{\ot}\mb{C}_1).\]
\end{definition}
\begin{theorem}[Kasparov]\label{T:ExtequalKK}
\[\EXT^{-1}(A,B) \cong KK^1(A,B).\]
\end{theorem}
\begin{proof}
 See \cite{kas80}. We note that an invertible extension $\tau$ in $\EXT^{-1}(A,B)$, which corresponds to a pair $(\phi,p)$, is mapped to a cycle $(H_B\oplus H_B, \phi \oplus \phi, 2p-1\oplus 1-2p)$ where $H_B\oplus H_B$ is graded by $1\oplus(-1)$ which we call the standard even grading.
\end{proof}

\begin{lemma}\label{L:com}
 Let $\phi:A \to B(H_{1}\oplus H_{2})$
 be a *-representation. Write $\phi(a) =
 \left(\begin{matrix}
 \phi_{11}(a)& \phi_{12}(a)  \\
 \phi_{21}(a) & \phi_{22}(a)
 \end{matrix} \right)$. Suppose $\phi_{11}$ is a *-homomorphism modulo
 $K(H_{1})$, i.e., $\dot{\phi}_{11}$ is a *-homomorphism. Then
  $\phi_{12}(a)$, $\phi_{21}(a)$ are compacts for any $a \in A$ and $\dot{\phi}_{22}$
 is a *-homomorphism.
 \end{lemma}
 \begin{proof}
  Using $\phi(aa^{*})=\phi(a)\phi(a^{*})$ under the decomposition of $\phi$ on
  $H_{1}\oplus H_{2}$ and the fact $\phi_{11}$ is *-homomorphism modulo
 $K(H_{1})$, we have $\phi_{12}(a)\phi_{12}(a^{*})$ is
 compact. Thus $\phi(a)$ is compact for any $a \in A$. Similarly, using
 $\phi(a^{*}a)=\phi(a)^{*}\phi(a)$, we have $\phi_{21}(a)$ is compact for any $a \in
 A$. It follows that $\phi_{22}$ is *-homomorphism modulo
 $K(H_2)$.
 \end{proof}
The following theorem due to Higson \cite{hi1} is important to us so that we give a complete proof.
\begin{proposition}\label{P:Ext }
 $K_{0}(D(A))\cong \EXT^{-1} ({A}, {\mathbb{C}})$.
 \end{proposition}
 \begin{proof}
 In light of Lemma \ref{L:K_0surjectivity},  we define the map from $K_{0}(D(A))$ to $\it
  {Ext}^{-1}(A,\mb{C})$ by  $$[p]_{0} \to [\tau_{\Phi,p}].$$ where $\Phi$ is an admissible representation of $A$
 on $H$ and $p \in D_{\Phi}(A)$.

  When $[p]_0=[q]_0$, as we have seen in the proof of Lemma \ref{L:K_0surjectivity}, $p$ and $q$
  are Murray-von Neumann equivalent in $D_{\Phi}(A)$ so that the partial isometry implementing
  this equivalence induces the equivalence between $\tau_{\Phi,p}$ and $\tau_{\Phi,q}$.
  Conversely, unitary equivalence between $\tau_{\Phi,p}$ and $\tau_{\Phi,q}$
  induces Murray-von Neumann equivalence between $p$ and $q$ evidently.\\
  From $\Phi\oplus\Phi \sim_{u} \Phi$, we get a unitary $u \in B({H\oplus H,H})$
  which induces  a natural isomorphism $Ad(u): \mathbb{M}_{2}(D_{\Phi}(A)) \to D_{\Phi}(A)$.
  Note that $\pi \circ Ad(u)=Ad(u)\circ (\pi\otimes id_{2})$.
  Since $[p]_0+[q]_0=[ p\oplus q]_0$ and  $p\oplus q  \in D_{\Phi\oplus\Phi}(A)$, $[p]_0+[q]_0$ is mapped to
$\left[ \pi \circ Ad(u) \circ  (\Phi\oplus\Phi)_{ p\oplus q}
\right]= \left[ Ad(u)\circ ((\pi\otimes id_{2})\circ
(\Phi\oplus\Phi)_{ p\oplus q}\right]$ which is indeed
$[\tau_{\Phi,p}]+[\tau_{\Phi,q}] $. So far we have shown the map is
a monomorphism.

 It remains to show the map is onto. Suppose $\rho$ is an invertible extension. Then $\rho$
is a semi-split extension of $A$ with a completely positive lifting
$\psi:A \to B(H)$. By the Stinespring's dilation
theorem, there is a non-degenerate *-representation $\phi: A \to
B(H_{0})$ and an isometry $V:H \to
H_{0}$ such that $\psi(a)=V^{*}\phi(a)V$ for all $a \in A$. It follows that $\ker (\dot{\phi})=0$ and $\phi$ is an admissible representation.
If we set $P_1=VV^{*}$ and $P_2=1-P_1$, then $H_{0}
=P_{1}(H_{0}) \oplus
P_2(H_{0})=H_{1}\oplus H_{2}$. If we
decompose $\phi$ on
$H_{0}=H_{1}\oplus H_{2}$ and write $\phi(a) =
 \left(\begin{matrix}
 \phi_{11}(a)& \phi_{12}(a)  \\
 \phi_{21}(a) & \phi_{22}(a)
 \end{matrix} \right)$, we have
$V\psi(a)V^{*}=VV^{*}\phi(a)VV^{*}=P_1\phi(a)P_1=\phi_{11}(a)$.
Since $\dot{\psi}=\rho$ is a (injective) *-homomorphism, we can
conclude $\phi_{11}$ is a injective *-homomorphism modulo compact.
By Lemma \ref{L:com}, $\phi_{12}(a)$, $\phi_{21}(a)$ are compacts
for $a \in A$ and $\phi_{22}$ is *a -homomorphism modulo compact.
This implies that $[P_1, \phi] \in K$. Thus
$\dot{\phi}_{11}$ is $\tau_{\phi, P_1}$. Viewing
$V:H \to P_1(H_{1})$ as a unitary, we can also see
that $\rho$ is unitarily equivalent to $\dot{\phi}_{11}$. Thus we
finish the proof.
\end{proof}

 \begin{theorem}[Odd case]\label{T:oddindex} The mapping $ K_{1}(B) \times K_{0}(D_{\Phi}(B)) \to
 \mathbb{Z} $ is the Kasparov product $ KK(S,B) \times KK^{1}(B,\mathbb{C}) \to
 \mathbb{Z}$.
\end{theorem}
\begin{proof}
Without loss of generality, we may assume that an element in $K_1(B)$ is represented by a unitary $v$ in B. Suppose an element in $K_0(D_{\Phi}(B))$ is represented by a projection $p$ in  $D_{\Phi}(B)$.

As we have noted, $[v]_1 \in K_1(B)$ is mapped to $\boldsymbol{\psi} \in KK(S,B)
$ where $\psi:S \to B$ is determined by sending $z-1$ to $v-1$. On
the other hand, $[p]$ is mapped to $[\tau_{\Phi,p}]$ by Proposition
\ref{P:Ext }. Using the isomorphism $\Ext {B} {\mathbb{C}} \to
KK^{1}(B, \mathbb{C})$ by Theorem \ref{T:ExtequalKK}, the image
of $\tau_{\Phi,p}$ is a cycle
$\mathcal{F}=\left(H\oplus H,\Phi\oplus\Phi,T\oplus-T\right)$
where $T= 2p-1$.

Then the Kasparov product $\boldsymbol{\psi}$ by
$[\mathcal{F}]$ is $$\boldsymbol{\psi}\cdot
[\mathcal{F}]=\left[\left(H\oplus H,(\Phi\circ\psi)\oplus(\Phi\circ\psi),T\oplus-T\right)\right
].$$
Note that  $(\Phi\circ\psi)(z-1)=\Phi(v-1)=\Phi(v)-1$. Under
the identification of $KK^{1}(S,\mathbb{C})$ with $
K_{1}(\mathcal{Q}(H))$, $\boldsymbol{\psi}\cdot
[\mathcal{F}]$ is mapped to $p\Phi(v)p-(1-p)$ by Proposition 17.5.7 in
\cite{bl}.  Using the index map
$\partial_1:K_{1}(\mathcal{Q}(H)) \to
K_0(K)=\mathbb{Z}$ in K-theory, we know
$\partial_1(p\Phi(v)p-(1-p))= \Ind (p\Phi(v)p-(1-p))$.
\end{proof}
 It is said that an element $\boldsymbol{x} \in \KK(A,B)$ is a
$\KK$-equivalence if there is a $\boldsymbol{y} \in \KK(B,A) $ such that
$\boldsymbol{x}\cdot \boldsymbol{y}=1_{A}$ and $\boldsymbol{y} \cdot
\boldsymbol{x}=1_{B}$. $A$ and $B$ are $KK$-equivalent if there is a
$\KK$-equivalence in $\KK(A,B)$. The following corollary, which is
originally due to Kasparov, establishes the $\KK$-equivalence of
$\mb{C}_1$ and $S$.

\begin{lemma}[Morita invariance]
Let $M_n$ denote the $n\times n$ matrix algebra of $\mb{C}$. Then
$$\KK(A,B)=\KK(A\otimes M_n, B\otimes M_m)$$
for any positive integers $n,m$.
\end{lemma}
\begin{lemma}
\[\KK(A\hat{\ot}\mb{C}_1, B\hat{\ot}\mb{C}_1)\cong KK(A,B)\]
\end{lemma}
\begin{proof}
Note that $\tau_{\mb{C}_1}\circ \tau_{\mb{C}_1}=\tau_{M_2}$. Thus Morita invariance implies that
$\tau_{\mb{C}_1}$ has the inverse as $\tau^{-1}_{M_2}\circ \tau_{\mb{C}_1}$.
\end{proof}
\begin{lemma}
\[ \KK(A\hat{\ot}\mathbb{C}_1,B)\cong \KK^1(A,B).\]
\end{lemma}
\begin{proof}
Recall that $KK^1(A,B)=KK(A, B \hat{\otimes}\mathbb{C}_1)$. Thus
\begin{align*}
 \KK(A \hat{\ot}\mathbb{C}_1,B) &\cong KK(A\otimes\mathbb{C}_1 \hat{\ot} \mb{C}_1, B\hat{\ot} \mb{C}_1)\\
                                &\cong KK(A \ot M_2, B\hat{\ot} \mb{C}_1)\\
                                &\cong KK(A,B\hat{\ot} \mb{C}_1 )\\
                                &\cong KK^1(A,B).
\end{align*}
\end{proof}
\begin{theorem}[Kasparov]
Let $\boldsymbol{x} \in \KK(\mathbb{C}_{1},S)\cong \KK^1(\mb{C},S)\cong
\EXT(\mathbb{C},S)$ be represented by the extension $ 0\to S \to
C \to \mathbb{C}\to 0 $ where $C$ is the cone and $\boldsymbol{y}
\in \KK(S,\mathbb{C}_{1})\cong \KK^1(S,\mb{C}) \cong \EXT(S,\mathbb{C})$ be represented
by the extension $0 \to K\to C^{*}(v-1)\to S \to 0 $ where
$v$ is a co-isometry of the Fredholm index 1 (e.g., the adjoint of
the unilateral shift on $H$). Then the Kasparov product
$\boldsymbol{x}\cdot \boldsymbol{y}=1_{\mathbb{C}_{1}}$.
\end{theorem}
\begin{proof}
  Note that $\boldsymbol{x}$ corresponds to the unitary $t \to e^{2\pi it}$
 in $K_{1}(S)$ by the Brown's Universal Coefficient Theorem \cite{br}. Also, the
 Busby invariant of $0 \to K\to C^{*}(v-1)\to S \to 0 $ is the
 extension $\tau: S \to \mathcal{Q}$ which sends $e^{2\pi it}-1$ to $\pi(v)-1$. Thus $\tau$ corresponds to a pair $(\Phi,1)$ such that $\Phi(e^{2 \pi i t})=v$.
If we apply Theorem \ref{T:oddindex} to the pair $([e^{2\pi
it}],[1]) \in K_1(S) \times K_0(D(S))$,  $\boldsymbol{x}\cdot \boldsymbol{y}$ is identified with $\Ind(v)=1$. Since $\KK(\mathbb{C}_{1},\mathbb{C}_{1})\cong \KK(\mb{C},\mb{C})\cong \mathbb{Z}$ by Proposition \ref{P:KKCC},
 we conclude $\boldsymbol{x}\cdot \boldsymbol{y}=1_{\mathbb{C}_{1}}$.
\end{proof}
Combining the above fact with a ring theoretic argument (See \cite{bl} for details), we can prove that $\boldsymbol{y} \cdot \boldsymbol{x}=1_{S}$.  Thus $\KK$-equivalence of $S$ induces the following corollary which will be used later.

\begin{corollary}\label{C:Bott}
$\tau_S: \KK(A,B) \to \KK(SA, SB)$ is an isomorphism.
\end{corollary}

  \begin{proposition}\label{P:KK}
 $ KK(A,\mathbb{C})\cong K_{1}(D_{\Phi})$ where $\Phi$ is an admissible representation of unital separable $C^{*}$-algebra $A$
 on a separable Hilbert space $H$.
 \end{proposition}
 \begin{proof}
 Using Lemma \ref{L:K_1surjectivity}, we define a map from $K_{1}(D_{\Phi})$ to $KK(A,\mathbb{C})$ by
 $$ [u]_{1} \to \left[\left( \hat{H},\left( \begin{matrix}\Phi & 0 \\
 0& \Phi \end{matrix}\right), \left( \begin{matrix}0 & u^{*} \\
 u& 0 \end{matrix}\right)  \right)\right]$$ where $\hat{H}$ is $H\oplus H$ with the standard
 even grading. Indeed, this construction gives
 rise to a well-defined group homomorphism. If $[u]=[v]$, then $u\oplus 1 $ is homotopic to $v\oplus 1$.

 Thus $\left( \hat{H}\oplus \hat{H} ,\left( \begin{matrix}\Phi\oplus\Phi & 0 \\ 0& \Phi\oplus\Phi \end{matrix}\right), \left( \begin{matrix}0 & u^{*}\oplus 1 \\
 u\oplus 1 & 0 \end{matrix}\right)  \right) \\ =
\left( \hat{H},\left( \begin{matrix}\Phi & 0 \\
 0& \Phi \end{matrix}\right), \left( \begin{matrix}0 & u^{*} \\
 u& 0 \end{matrix}\right)  \right)\oplus \left( \hat{H},\left( \begin{matrix}\Phi & 0 \\
 0& \Phi \end{matrix}\right), \left( \begin{matrix}0 & 1 \\
 1 & 0 \end{matrix}\right)  \right)$ is operator homotopic to $\left( \hat{H},\left( \begin{matrix}\Phi & 0 \\
 0& \Phi \end{matrix}\right), \left( \begin{matrix}0 & v^{*} \\
 v& 0 \end{matrix}\right)  \right)\oplus \left( \hat{H},\left( \begin{matrix}\Phi & 0 \\
 0& \Phi \end{matrix}\right), \left( \begin{matrix}0 & 1 \\
 1 & 0 \end{matrix}\right)  \right) \\ = \left( \hat{H}\oplus \hat{H},\left( \begin{matrix}\Phi\oplus\Phi & 0 \\ 0& \Phi\oplus\Phi \end{matrix}\right), \left( \begin{matrix}0 & v^{*}\oplus 1 \\
 v\oplus 1 & 0 \end{matrix}\right) \right)$. Similarly, it can be shown that the map is a group homomorphism.
 If we show that it is surjective, we are done. We will use Higson's idea in p354 \cite{hi1}.
 Let $\alpha \in KK(A,\mathbb{C})$ be represented by
 $\left(H_{0}\oplus H_{1},\left( \begin{matrix}
 \phi_{0}&0 \\
 0 & \phi_{1}
 \end{matrix}\right),
 \left( \begin{matrix}
 0 & u^{*} \\
 u & 0
 \end{matrix}\right) \right)$ where $u$ is a unitary in $B(H_{0},H_{1})$.
 Let $\Psi= \cdots \oplus \phi_{0}\oplus\phi_{0}\oplus \phi_{1}\oplus\phi_{1}\oplus \cdots$ and
 ${\bf H}= \cdots \oplus H_{0}\oplus H_{0}\oplus H_{1}\oplus H_{1}\cdots \oplus$.
 We consider a degenerate cycle $$\left({\bf H}\oplus {\bf H},\left(
 \begin{matrix}
 \Psi& 0 \\
 0 & \Psi
 \end{matrix}\right),\left(
 \begin{matrix}
 0 & I \\
 I & 0 \end{matrix}\right) \right) $$ Then
 $$\left(H_{0}\oplus H_{1},\left(
 \begin{matrix}
 \phi_{0}&0\\
 0 & \phi_{1}
 \end{matrix}\right),\left(
 \begin{matrix}
 0 & u^{*} \\
 u & 0
 \end{matrix}\right) \right)\oplus \left({\bf H}\oplus {\bf H},\left(
 \begin{matrix}
 \Psi& 0 \\
 0 & \Psi \end{matrix}\right),\left(
 \begin{matrix}
 0 & I \\
 I & 0 \end{matrix}\right) \right)$$ is unitarily equivalent to $\left({\bf H}\oplus {\bf H},
 \left( \begin{matrix}
 \Psi& 0 \\
 0 & \Psi
 \end{matrix}\right),\left(
 \begin{matrix}
 0 & F^{*} \\
 F & 0
 \end{matrix}\right) \right)$ where $F= \left( \begin{matrix}I & 0& 0\\ 0&u&0\\ 0&0& I \end{matrix}\right)\circ $
 shifting to the right, i.e., $F$ sends $(\cdots,\eta_{1},\eta_{0},\xi_{0},\xi_{1},\cdots)$
 to $(\cdots,\eta_{2},\eta_{1},u\eta_{0},\xi_{0},\cdots)$.\\
 Again by adding a degenerate cycle $\left( \hat{H},\left( \begin{matrix}\Phi & 0 \\
 0& \Phi \end{matrix}\right), \left( \begin{matrix}0 & I \\
 I & 0 \end{matrix}\right) \right)$, we get
 \begin{align*}\alpha &=\left[\left({\bf H}\oplus H,\left(
 \begin{matrix}
 \Psi\oplus \Phi & 0 \\
 0& \Psi\oplus\Phi
 \end{matrix}\right),
 \left( \begin{matrix}
 0 & F^{*}\oplus I \\
 F\oplus I & 0
 \end{matrix}\right)  \right)\right]\\
  \intertext{ Since $\Phi$ is admissible, we can have a unitary $U \in
  B({\bf H}\oplus H,H)$ such that $\Omega=Ad(U)\circ
  (\Psi\oplus\Phi) \sim \Phi $.}
&=\left[ \left(\hat{H},\Omega\oplus\Omega, \left(
\begin{matrix}
0 & Ad(U)(F^{*}\oplus I) \\
 Ad(U)(F\oplus I) & 0
 \end{matrix}\right) \right)\right]\\
 \intertext{By Lemma 4.1.10. in \cite{jeth}}
 &=\left[ \left(\hat{H},\Phi\oplus \Phi,
 \left( \begin{matrix}
 0 & Ad(U)(F^{*}\oplus I) \\
 Ad(U)(F\oplus I) & 0
 \end{matrix}\right) \right)\right]
 \end{align*}
 It is not hard to check that  $Ad(U)(F\oplus I)\in D_{\Phi}(A)$ and $\alpha$ is the image of it.
 Thus we finish the proof.
 \end{proof}
 \begin{theorem}[Even case]\label{T:evenindex} The mapping $ K_{0}(A) \times K_{1}(D_{\Phi}) \to
 \mathbb{Z} $ is the Kasparov product $KK(S,SA) \times KK(SA,S) \to
\mathbb{Z} $.
 \end{theorem}
 \begin{proof}
 Without loss of generality, we
 may assume $p$ is the element of $A$. (If necessary, consider
 $\mathbb{C}^{k}\otimes A $.) Using the Bott map in K-theory, $p$
 gets mapped to $f_{p}(z) \in K_{1}(SA)$. Then, as we have noted in Proposition \ref{P:1}, $f_{p}(z)$ is lifted to
 $\boldsymbol{\Psi}$ as an element of $KK(S,SA)$ where $\Psi$ is the *-homomorphism from $S$ to $SA$ that is determined by sending $z-1$ to $(z-1)p$.\\
 On the other hand, $[u] \in K_{1}(D_{\Phi}(A))$ is mapped to $[\mathcal{E}]$ in $KK(A,\mathbb{C})$ where
 $\mathcal{E}=\left( \hat{H},\left( \begin{matrix}\Phi & 0 \\
 0& \Phi \end{matrix}\right), \left( \begin{matrix}0 & u^{*} \\
 u& 0 \end{matrix}\right)  \right)$ by
 Proposition \ref{P:KK}.\\
 Using natural isomorphism $\tau_{S}:KK(A,\mathbb{C})\to KK(SA,S)$, we can think of a Kasparov product $\boldsymbol{\Psi}$ by $[\tau_{S}(\mathcal{E})]$.\\
 Using  elementary functorial properties, we can check  $\boldsymbol{\Psi}\cdot[\tau_{S}(\mathcal{E})]$ is equal to $\left[(\hat{H}\otimes S,((\Phi\oplus\Phi)\otimes 1) \circ \psi, G\otimes 1 ) \right]$ denoting  $\left( \begin{matrix}0 & u^{*} \\
 u& 0 \end{matrix}\right)$ by $G$.\\
Since $\tau_{S}:KK(\mathbb{C},\mathbb{C})\to KK(S,S)$ is an isomorphism, there
is a map $\rho:\mathbb{C}\to B(H)$ such that
\begin{equation*}\label{E:first}
\tau_{S}((\hat{H},\rho\oplus\rho, G ))=
(\hat{H}\otimes S,((\Phi\oplus\Phi)\otimes 1) \circ \psi,
G\otimes 1 ).
\end{equation*}
This implies that $(\rho \oplus \rho) \ot 1 = ((\Phi\oplus\Phi)\otimes 1) \circ \psi$. Thus
\begin{align*}
(\rho\otimes 1) (z-1)&=(z-1)\rho(1)\\
 &=(\Phi\otimes 1) (\Psi(z-1)) \quad \text{by equality (\ref{E:first})}\\
 &=(\Phi\otimes 1) ((z-1)p)\\
 &=(z-1)\Phi(p).
\end{align*}
Consequently, $\rho(1)=\Phi(p)$.\\
Then $\left(\hat{H},\rho\oplus\rho, G \right)$ is mapped to
$\Ind(\Phi(p)u\Phi(p))$ by Proposition \ref{P:KKCC}.
 \end{proof}


\begin{thebibliography}{99}
 \bibitem[Ar]{ar} W. Arveson \emph{Notes on extensions of $C^{*}$-algebras}
 Duke Math. Journal Vol44 no2 329-355 (1977)
 \bibitem[Br]{br} L. G. Brown \emph{The universal coefficient theorem for Ext and quasidiagonality} pp60-64 in Operator algebras and group representations, I. Monographs of Stud. Math. 17, Pitman, Boston, Mass. 1984
 \bibitem[Bl]{bl} B. Blackadar \emph{K-theory for Operator
 Algebras} MSRI Publicatons Vol.5 Second Edition Cambridge University Press
 1998
 \bibitem[Hig]{hi1} N. Higson  \emph{$C^{*}$-Algebra Extension Theory and Duality}
       J. Func. Anal. 129(1995), 349-363.
 \bibitem[JenThom]{jeth} K. N. Jensen, K. Thomsen \emph{Elements of KK-
 theory} Birkh\"{a}user, Boston, 1991 MR94b:19008
 \bibitem[Kas]{kas80} G. G. Kasparov \emph{The Operator K
 functor and extensions of $C^{*}$-algebras} Mathe. USSR-Izv
 16(1981) 513-572[English Tanslation]

 \bibitem[Pa]{pa} W. Pascheke \emph{K-theory for commutants in the Calkin algebra}
 Pacific J. Math. 95 (1981) 427-437
 \end{thebibliography}
\end{document}